\documentclass[12pt,reqno]{amsart} 
\usepackage{amsfonts,amsmath,amssymb}
\usepackage[english]{babel}  
\usepackage[all]{xy}

\swapnumbers

\newcommand{\Z}{{\mathbb Z}}

\renewcommand{\P}{{\mathbb P}}
\renewcommand{\H}{{\text{H}}}

\newcommand{\ke}{{\mathcal E}}
\newcommand{\kf}{{\mathcal F}}

\newcommand{\ki}{{\mathcal I}}

\newcommand{\ko}{{\mathcal O}}

\DeclareMathOperator{\Ext}{Ext}

\begin{document}

\title[Babylonian Tower]{The Splitting Criterion of Kempf and the\\
 Babylonian Tower Theorem}

\author{I.~Coand\u{a}} \address{
  Institute of Mathematics
  of the Romanian Academy\\
  P.O.~Box 1-764
\newline RO--70700 Bucharest, Romania}
\email{Iustin.Coanda@imar.ro} 
\author{G.~Trautmann}
\address{
  Universit\"at Kaiserslautern\\
  Fachbereich Mathematik\\
  Erwin-Schr\"odinger-Stra{\ss}e
\newline D-67663 Kaiserslautern}
\email{trm@mathematik.uni-kl.de}

\footnotetext{Research supported by the DFG-Schwerpunktprogramm 1094}

\begin{abstract}

\vskip5mm
We show that the idea used by Kempf (1990)  in order to obtain a
splitting criterion for vector bundles on projective spaces leads to an
elementary proof of the Babylonian tower theorem for this class of bundles, a
result due to Barth--Van de Ven (1974) in the rank 2 case and to 
Sato (1977) and Tyurin (1976) in the case of arbitrary rank. 
As a byproduct we obtain a slight improvement of the numerical
criterion of Flenner (1985) in the particular case under consideration. 

\end{abstract}

\maketitle

\vskip1cm 
The aim of this note is to provide an elementary and short proof of the
following Babylonian tower theorem.

{\bf Theorem:}\; {\em Let $\ke$ be a locally free sheaf on the projective space
  $\P_n$ (defined over an algebraically closed field $k$), $n\geq 2$. If
$m>\Sigma_{i>0} \dim \Ext^1(\ke, \ke (-i)),$
then $\ke$ cannot be extended to $\P_{n+m}$ unless it is a direct sum of 
line bundles.}

The Babylonian tower theorem for vector bundles on a punctured spectrum of
H. Flenner, \cite{fl}, asserts, in our special case, the same conclusion but 
for the larger bound
$$
m>\Sigma_{i\in\Z} \dim \Ext^1(\ke, \ke(-i)).
$$
On the other hand, Kempf, \cite{ke}, proved that $\ke$ cannot be extended to
$\P_{n+1}$ if $\Ext^1(\ke, \ke(-i))=0$ for all $i>0$, unless it splits into a
direct sum of line bundles. Our proof of the theorem is inspired by Kempf's
proof leading us to use infinitesimal neighbourhoods. The idea of Kempf is
reflected in the following

{\bf Lemma:}\; {\em Let $\kf$ be locally free on $\P_n, n\geq 2$, let 
$H\subset \P_n$
be a hyperplane of equation $h=0$, and let $H_i$ be the $i$--th
  infinitesimal neighbourhood of equation $h^{i+1}=0$. The linear
  projection onto $H$ from a point of $\P_n\smallsetminus H$ endows
  $\ko_{H_i}$ with the structure of an $\ko_H$--algebra. Then, if for every
  $i\geq 1$ the exact sequence
$$
o\to \kf_{H_{i-1}} (-1)\xrightarrow{h} \kf_{H_i}\to \kf_H\to 0
$$ 
splits as a sequence of $\ko_H$--modules, $\kf$ is a direct sum of line 
bundles.}

\begin{proof} Here and elsewhere $\kf_Y$ shall denote the restriction of $\kf$
to a subscheme $Y$.
A splitting of the sequence of $\kf_{H_i}$ implies that the induced map
$\H^p\kf_{H_i} (a)\to \H^p\kf_H(a)$
is surjective for every $p\geq 0,\; i\geq 0,\; a\in\Z$. Using the exact sequence
$$
0\to \kf(-i-1)\xrightarrow{h^{i+1}} \kf\to\kf_{H_i}\to 0
$$
one deduces that for $p\leq n-2$ and any $a\in\Z$ the map
$\H^p\kf(a)\to \H^p\kf_{H_i}(a)$
is surjective for large $i>>0$, because $\H^{p+1}\kf(a-i-1)=0$ for large
$i$. Consequently, also the maps
$\H^p\kf(a)\to \H^p\kf_H(a)$
are surjective for any $p\leq n-2$ and any $a\in\Z$.
It follows that the maps $\H^p\kf(a-1)\xrightarrow{h} \H^p\kf(a)$ are injective
for $1\leq p\leq n-1$ and any $a\in\Z$. But $\H^p\kf(a)=0$ for $a>>0$. Then
$\H^p\kf(a)=0$ for every $1\leq p\leq n-1$ and $a\in\Z$. By the well known
splitting criterion of Horrocks, $\kf$ is a direct sum of line bundles.
\end{proof}
\vskip5mm

Before proving the theorem we do some preparations.
Let $\P_n$ be embedded in $\P=\P_{n+m}$ as the linear subspace $L$ with 
equations
$x_{n+1}=\cdots=x_{n+m}=0$ and let $S=k[x_0, \ldots, x_{n+m}]$ be the
homogeneous coordinate ring of $\P_{n+m}$. Let $\ki$ denote the ideal sheaf of
$L$ and let $L_i$ be the $i$--th
infinitesimal neighbourhood of $L$ in $\P_{n+m}$ defined by $\ki^{i+1}$. Let
$\P_{n+m}\smallsetminus L'\xrightarrow{\pi}L$ be the central projection, where
$L'=\{x_0=\cdots = x_n=0\}$ is complementary to $L$ with homogeneous
coordinate ring $R:=k[x_{n+1}, \ldots, x_{n+m}]$. Then $\pi$ endows $\ko_{L_i}$
with the structure of an $\ko_L$--algebra. Now, as an $\ko_L$--module, 
$$
\ko_{L_i}\cong \ko_L\oplus \ko_L(-1)\otimes_k R_1
\oplus\cdots\oplus \ko_L(-i)\otimes_k R_i 
$$
where $R_i$ are the graded components of $R$, noting that $\pi$ corresponds to the 
inclusion $k[x_0, \ldots, x_n]\hookrightarrow S$. Moreover, the ring structure of 
$\ko_{L_i}$ corresponds to the natural
multiplication in this decomposition, defined by
$$
(f_p\otimes r_p)(f_q\otimes r_q)=f_pf_q\otimes r_pr_q.
$$
This allows to consider the following ideals of $\ko_{L_i}$. If $J$ is a
homogeneous ideal of $R$, then 
$$J^{(i)}:= \ko_L(-1)\otimes_k
  J_1\oplus\ldots\oplus\ko_L(-i)\otimes_k J_i
$$
is an ideal subsheaf of $\ko_{L_i}$. In fact, if $X\subset \P_{n+m}$ is the
subscheme defined by the ideal $JS$, then
$\ko_{L_i}/J^{(i)}\cong \ko_{L_i\cap X}.$
In particular, if $J$ is the ideal of a point $p\in L'$, then $JS$ defines the
linear subspace $P\subset \P_{n+m}$, which is the linear span of $L$ and $p$, 
and in this case $\ko_{L_i}/J^{(i)}\cong \ko_{L_i\cap P}$.
\vskip1cm

{\em Proof of the theorem.} 
Let $\ke$ be as in the theorem and suppose that there exists a locally free
sheaf $\kf$ on $\P_{n+m}$ such that $\kf_L\cong\ke$. We shall construct a
homogeneous ideal $J\subset R$, generated by (at most)
$\Sigma_{i>0}\dim\Ext^1(\ke, \ke(-i))$ homogeneous elements, such that
for any $i\geq 0$ the short exact sequence 
$$
0\to \kf_{L_i}\otimes_{\ko_{L_i}}\ki\ko_{L_i}/ J^{(i)}\to \kf_{L_i}\otimes_{\ko_{L_i}}
\ko_{L_i}/J^{(i)}\to\kf_L\to 0 \eqno (A_i)
$$
splits as a sequence of $\ko_L$--modules. Since $m-1\geq
\Sigma_{i>0}\dim\Ext^1(\ke, \ke(-i))$, there exists a point $p\in L'\cong \P_{m-1}$ such
that the elements of $J$ vanish in $p$. Let $P$ be the span of $p$ and $L$ as
in the preparation. We obtain the commutative diagram
$$
\xymatrix{0\ar[r] & \kf_{L_i}\otimes_{\ko_{L_i}} \ki\ko_{L_i}/J^{(i)}\ar[d]\ar[r] & 
\kf_{L_i}\otimes_{\ko_{L_i}} \ko_{L_i}/J^{(i)} \ar[d]\ar[r] & \kf_L\ar@{=}[d]\ar[r] & 0\\
0\ar[r] & \kf_{L_{i-1}\cap P}(-1)\ar[r] & \kf_{L_i\cap P}\ar[r] & \kf_L\ar[r] & 0}
$$
which implies that the bottom row splits, too. Because the schemes $L_i\cap P$ 
are the infinitesimal neighbourhoods
of $L$ in $P$, the Lemma implies that $\kf_P$ is a direct sum of line bundles,
and so is $\ke\cong\kf_L$.
\vskip3mm

As for $J$, it is obtained as the union of an ascending chain of homogeneous
ideals of $R$, constructed inductively by a standard technique from
infinitesimal deformation theory, adapted to our situation. The induction step
is as follows: Let $i\geq 0$ be an integer and assume that an ideal 
$J\subset R$ has been constructed such that the sequence $(A_j)$ splits for 
$j\leq i$. Consider the naturally induced diagram
$$
\xymatrix{0\ar[r] & \kf_{L_{i+1}}\otimes_{\ko_{L_{i+1}}} \ki\ko_{L_{i+1}}/J^{(i+1)}\ar[d]\ar[r] &
\kf_{L_{i+1}}\otimes_{\ko_{L_{i+1}}} \ko_{L_{i+1}}/J^{(i+1)} \ar[d]\ar[r] &
\kf_L\ar@{=}[d]\ar[r] & 0\\
0\ar[r] & \kf_{L_i}\otimes_{\ko_{L_i}} \ki\ko_{L_i}/J^{(i)}\ar[r] &
\kf_{L_i}\otimes_{\ko_{L_i}} \ko_{L_i}/J^{(i)} \ar[r] &
\kf_L\ar[r] & 0}
$$
and let $e\in\Ext^1_{\ko_L}(\kf_L, \kf_{L_{i+1}}\otimes\ki\ko_{L_{i+1}}/J^{(i+1)})$ 
be the extension class (over $\ko_L$) of the top row, which is
then mapped in $\Ext^1_{\ko_L}(\kf_L, \kf_{L_i}\otimes\ki\ko_{L_i}/J^{(i)})$ to 
the extension class of the bottom row, which is supposed to be
$0$. Since we have the exact sequence
$$
0\to\ko_L(-i-1)\otimes_k R_{i+1}/J_{i+1}\to
\ki\ko_{L_{i+1}}/J^{(i+1)}\to \ki\ko_{L_i}/J^{(i)}\to 0,
$$
it follows that $e$ is the image of an element
$e'\in \Ext^1_{\ko_L}(\kf_L, \kf_L(-i-1))\otimes_k R_{i+1}/J_{i+1}.$
Let $e_1, \ldots, e_s$ be a $k$--basis of $\Ext^1_{\ko_L}(\kf_L,
\kf_L(-i-1))$. Then $e'=e_1\otimes \bar{f}_1+\ldots+e_s\otimes \bar{f}_s$ with
$f_1, \ldots, f_s\in R_{i+1}$. If $K:=J+Rf_1+\cdots +Rf_s$, then $K_j= J_j$ for
$j\leq i$, such that the sequences $(A_j)$ for $K$ split for $j\leq i$. In
addition,
$$
0\to \kf_{L_{i+1}}\otimes \ki\ko_{L_{i+1}}/K^{(i+1)}\to \kf_{L_{i+1}}\otimes
  \ko_{L_{i+1}}/K^{(i+1)}\to\kf_L\to 0
$$
splits, too, because of the commutativity of the diagram
$$
\xymatrix{
\Ext^1_{\ko_L}(\kf_L, \kf_L(-i-1))\otimes_k R_{i+1}/J_{i+1}\ar[r]\ar[d] &
\Ext^1_{\ko_L}(\kf_L, \kf_L(-i-1))\otimes_k R_{i+1}/K_{i+1}\ar[d]\\
\Ext^1_{\ko_L}(\kf_L, \kf_{L_{i+1}}\otimes\ki\ko_{L_{i+1}}/J^{(i+1)})\ar[r] &
\Ext^1_{\ko_L}(\kf_L, \kf_{L_{i+1}}\otimes\ki\ko_{L_{i+1}}/K^{(i+1)})
}
$$
and of the fact that $e'$ is mapped to $0$ in 
$\Ext^1_{\ko_L}(\kf_L, \kf_L(-i-1))\otimes_k R_{i+1}/K_{i+1}.$
This completes the proof of the theorem.

\end{document}